\pdfoutput=1

\documentclass[12pt,reqno]{amsart}

\usepackage{amsmath,amsthm,amssymb}

\usepackage{graphicx}

\newcommand{\delthe}{4/(\delta-\theta)}

\newtheorem{theorem}{Theorem}[section]

\newtheorem{corollary}[theorem]{Corollary}

\newtheorem{definition}[theorem]{Definition}

\newtheorem{lemma}[theorem]{Lemma}

\newtheorem{proposition}[theorem]{Proposition}

\numberwithin{equation}{section}

\theoremstyle{definition}
\newtheorem{remark}[theorem]{Remark}

\newcommand\ben{\begin{enumerate}}
\newcommand\een{\end{enumerate}}

\renewcommand{\i}{{\mathrm{i}}} 

\newcommand{\twocase}[5]{#1 \begin{cases} #2 & \text{#3}\\ #4
&\text{#5} \end{cases}   }

\newcommand{\threecase}[6]{
\begin{cases} 
#1 & \text{#2}\\ 
#3 & \text{#4}\\
#5 & \text{#6}\\ 
\end{cases}   }

\newcommand\be{\begin{equation}}
\newcommand\ee{\end{equation}}
\newcommand\benn{\begin{equation*}}
\newcommand\eenn{\end{equation*}}
\newcommand\bea{\begin{eqnarray}}
\newcommand\eea{\end{eqnarray}}
\newcommand\beann{\begin{eqnarray*}}
\newcommand\eeann{\end{eqnarray*}}

\renewcommand{\H}{\mathbb{H}}
\newcommand{\boldH}{\mathbb{H}}
\newcommand{\R}{\mathbb{R}}

\newcommand{\Z}{\mathbb{Z}}

\newcommand{\N}{\mathbb{N}}
\newcommand{\F}{\mathbb{F}}

\newcommand{\B}{\frak{B}}

\newcommand{\Fu}{\mathcal{F}} 
\newcommand{\A}{\mathcal{A}}
\newcommand{\X}{\mathcal{X}}

\newcommand{\Or}{\ensuremath{{\mathcal O}}}

\newcommand{\ga}{\alpha}     
\newcommand{\gb}{\beta}      
\newcommand{\gd}{\delta}     
\newcommand{\gD}{\Delta}     
\newcommand{\gep}{\epsilon}  
\newcommand{\vep}{\varepsilon} 
\newcommand{\G}{\Gamma}      
\newcommand{\Gi}{\Gamma_{\infty}}      
\newcommand{\g}{\gamma}      
\newcommand{\gL}{\Lambda}    
\newcommand{\gl}{\lambda}    
\newcommand{\gk}{\kappa}    
\newcommand{\gS}{\Sigma}     
\newcommand{\gt}{\theta}     
\newcommand{\gw}{\omega}      

\newcommand{\foh}{\frac{1}{2}}  

\newcommand{\vectwo}[2]
{\left(\begin{array}{c}
                        #1    \\
                        #2
                          \end{array}\right) }

\newcommand{\mattwo}[4]
{\left(\begin{array}{cc}
                        #1  & #2   \\
                        #3 &  #4
                          \end{array}\right) }

\newcommand{\Spec}{\operatorname{Spec}}

\newcommand{\PSL}{\operatorname{PSL}}
\newcommand{\SL}{\mbox{SL}}

\newcommand{\vol}{\operatorname{vol}}
\newcommand{\tr}{\mbox{tr}}
\renewcommand{\mod}{\operatorname{mod}}

\renewcommand{\hat}{\widehat} 

\renewcommand{\Re}{{\mathfrak{Re}}}
\renewcommand{\Im}{{\mathfrak{Im}}}

\newcommand{\<}{\left\langle}
\renewcommand{\>}{\right\rangle}
\newcommand{\bk}{\backslash}
\newcommand{\GbkH}{\G\bk\boldH}
\newcommand{\GibkH}{\Gi\bk\boldH}

\newcommand{\ddxt}{{\partial^2\over\partial x^2}}
\newcommand{\ddyt}{{\partial^2\over\partial y^2}}

\begin{document}

\title[Lattice Point Count and Sieves]{The Hyperbolic Lattice Point Count  in Infinite Volume with Applications to Sieves}
\author{Alex V. Kontorovich} \email{alexk@math.brown.edu}
\begin{abstract}

We develop novel techniques using 
 abstract operator theory to obtain asymptotic formulae for lattice counting problems on infinite-volume hyperbolic manifolds, with error terms which are uniform as the lattice moves through ``congruence'' subgroups. 
We give the following application to the theory of affine linear sieves. 
 In the spirit of Fermat, consider the problem of primes in the sum of two squares, $f(c,d)=c^2+d^2$, but restrict $(c,d)$ to the orbit $\Or=(0,1)\G$, where $\G$ is an infinite-index non-elementary finitely-generated subgroup of $\SL(2,\Z)$. Assume that the Reimann surface $\G\bk\boldH$ has a cusp at infinity.  
  We show that the set of values $f(\Or)$ contains infinitely many integers having at most $R$ prime factors for any $R>\delthe$, where $\gt>1/2$ is the spectral gap and $\gd<1$ is the Hausdorff dimension of the limit set of $\G$. If $\gd>149/150$, then we can take $\gt=5/6$, giving $R=25$. The limit of this method is $R=9$ for $\gd-\gt>4/9$. This is the same number of prime factors as attained in Brun's original attack on the twin prime conjecture.
\end{abstract}
\address{Department of Mathematics,
Columbia University, 2990 Broadway, New York, NY $10027$} 
\curraddr{Department of Mathematics,
Brown University, 151 Thayer St, Providence, RI $02912$} 

\subjclass[2000]{Primary: 11N32, 30F35; Secondary: 11F72, 11N36}

\keywords{Affine Linear Sieve, Fuchsian Groups, Lattice Point Count, Patterson-Sullivan Theory} 

\date{\today}
\maketitle
\tableofcontents

\section{Introduction}

Many of the most enticing problems in number theory amount to finding primes or almost-primes (numbers having few prime factors) in ``thin'' subsets of the integers. By thin, we mean that the number of elements in the given set of size not exceeding a height $T$ is bounded by $T^{\gd}$, where $\gd$ is some constant less than $1$. See the papers \cite{PiatetskiShapiro1953, Chen1973, Iwaniec1978, FriedlanderIwaniec1998, HeathBrown2001} which are landmarks in sieve theory, producing primes
or products of at most two primes in thin sets.

It is our main goal to exhibit such a set arising from the orbit of an affine linear group action and sift 
for elements having few prime factors. Our starting point is the fundamental work of Bourgain, Gamburd, and Sarnak \cite{BourgainGamburdSarnak2006, BourgainGamburdSarnak2008, Sarnak2007} on the so-called  affine linear sieve:

\begin{theorem}[\cite{BourgainGamburdSarnak2006}]\label{BGS}
Let $\G\subset\SL_2(\Z)$ be any non-elementary\footnote
{Recall that an elementary group $\G\subset\SL_2(\Z)$ has the property that any two elements $g,h\in\G$ of infinite order have $|\tr(ghg^{-1}h^{-1})|=2$, see \cite{Beardon1983}. These are the ``abelian'' groups, such as purely hyperbolic groups, or purely elliptic groups (a torus generated by a single elliptic element).} 
subgroup, let $q\in\Z^2$ be nonzero, and let $\Or=q\cdot\G$ 
be a $\G$-orbit. Let $f:\Z^2\to\Z$ be any polynomial. Then there exists an $R<\infty$, depending on all of the above data, such that there are infinitely many points in the set  $f(\Or)$ having at most  $R$ prime factors. 
\end{theorem}

Notice that there are no congruence conditions, nor conditions on $f$ (such as irreducibility) -- these are all factored into $R$, which is left completely unspecified. There are three main 
ingredients:
\begin{enumerate}
\item Number Theory -- the combinatorial sieve,
\item Algebra -- Strong Approximation and Goursat's Lemma, and
\item Combinatorics -- counting by wordlength in $\G$ 
and 
extending the recently established
expander property of Bourgain-Gamburd \cite{BourgainGamburd2007} 
(with an unspecified spectral gap)
to square-free moduli.
\end{enumerate}

The lack of specificity of the spectral gap (and arbitrariness of the choice of orbit $\Or$ and function $f$) induces a lack of specificity of the number $R$ of prime factors.

In this paper, we select a particular orbit $\Or$ and function $f$, and give
 a precise bound for the number  $R$ of prime factors  by replacing input (3) above   with an archimedian count and using Gamburd's explicit $5/6$-th gap \cite{Gamburd2002} in place of the unspecified spectral gap. 
 
 This amounts  to a hyperbolic lattice point counting problem, requiring uniform error estimates as the lattice moves through ``congruence'' subgroups. In the interest of having a thin set, we are compelled to work in infinite volume, where standard spectral methods (decomposition into Maass forms and Eisenstein series) cannot be applied.
  Therefore we  develop novel ``soft'' methods using only  operator theory (the abstract spectral theorem; see e.g. \cite{Halmos1963}) and spectral information from Patterson-Sullivan theory \cite{Patterson1976,Sullivan1984} and Lax-Phillips \cite{LaxPhillips1982} to circumvent explicit knowledge of a spectral decomposition.

\begin{remark}
 Note that
an archimedean count is also used in Theorem 2 of
\cite{BourgainGamburdSarnak2006}; Corollary \ref{corcor} in the present paper is closely related to
this result. The difference is counting in a group versus counting in an orbit.
 The problem of counting in the group in infinite volume was solved in Lax-Phillips \cite{LaxPhillips1982}. 
To count in an orbit, one faces the serious issue of a stabilizer, as we discuss below.
\end{remark}

To state our main theorem, we require some notation. Our function of choice will be the sum of two squares: 
\be\label{fIs}
f(c,d):=c^2+d^2. 
\ee
Let $\G\subset\SL_2(\Z)$ be a non-elementary finitely-generated Fuchsian group. We will soon turn our attention exclusively to groups of the second kind -- ones having infinite co-volume -- but we do not make this restriction just yet. 
Denote by $\Gi$ the set of 
elements of $\G$ which stabilize infinity:
\be\label{def:Gi}
\Gi:=\left\{\g\in\G\mid\g=\mattwo{1}{*}{0}{1}\right\}.
\ee
Let $\gd=\gd(\G)
\le1
$ denote the Hausdorff dimension of the limit set of $\G$ and $\gt
<
\gd
$ be the spectral gap
 (see \S \ref{delta}). The index of $\G$ in $\SL_2(\Z)$ is finite if and only if $\gd=1$. 
 
 Let $\Or$ be the orbit of bottom rows of $\G$, 
\be\label{def:O}
\Or:=(0,1)\G = \left\{(c,d):\mattwo{*}{*}{c}{d}\in\Gi\bk\G\right\},
\ee
and for  a height $T>1$  let $\Or(T)$ denote the set of orbital points not exceeding this height:
\be\label{def:OT}
\Or(T):=\{(c,d)\in\Or\mid c^2+d^2<T\}.
\ee
For $R\ge 1$ let $\Or(T,R)$ be defined by
\benn
\Or(T,R):=\{(c,d)\in\Or(T) \mid f(c,d) \text{ has at most $R$ prime factors}\}.
\eenn
Recall the notation 
$$
f\asymp g\quad\text{ for \quad $g\ll f \ll g$.}
$$

\begin{theorem}[Main Theorem] \label{mainThm}
Let $\G\subset \SL_2(\Z)$ be a non-elementary finitely-generated Fuchsian group, let $\gd$ be the Hausdorff dimension of its limit set, $\gt$  the spectral gap, and $\Or(T)$ and $\Or(T,R)$ defined as above. Assume  $\Gi$  is nontrivial (then $\gd>1/2$, \cite{Beardon1968}). 
As $T\to\infty$,
\begin{enumerate}
\item There exist  constants $c_0>0$  and $\eta>0$ such that 
\be\label{mainthm1}
|\Or(T)| =  c_0T^{\gd} + O(T^{\gd-\eta}).
\ee
\item For any fixed $R>4/(\gd-\gt)$,
\be\label{Part2}
 |\Or(T,R)|\asymp T^{\gd}/\log T.
\ee
\end{enumerate}
\end{theorem}

Part (1) above
tells us that the orbit is
 thin 
 if and only if $\gd<1$, and is the quintessence of our necessity to work in infinite co-volume. 
In \cite{BourgainGamburdSarnak2006}, it is proved that if $\gd>1/2$, then there always exists {\it some} spectral gap $\gt<\gd$. Gamburd \cite{Gamburd2002} shows that if $\gd>5/6$ then we can take $\gt=5/6$. If
on the other hand
 $\G$
is a finite co-volume group and moreover
 a congruence group, then Kim-Sarnak \cite{KimSarnak2003} allows 
 $\gt=1/2+7/64
\approx
.609
$. The following corollary is immediate.

\begin{corollary}
As $T\to\infty$,
\benn
|\Or(T,R)|\gg T^{\gd}/\log T
\eenn
 for
\begin{enumerate}
\renewcommand{\labelenumi}{(\roman{enumi})}
\item $R=25$ if $\gd>149/150$ (setting $\gt=5/6$), 
\item $R=11$ if $\G$ is a finite index congruence group  (where $\gd=1$  and    $\gt =39/64$), and 
\item $R=9$ if $\gd-\gt>4/9$.
\end{enumerate}
\end{corollary}

It is known that in dimension two, infinite co-volume groups exist with $\gd$ arbitrarily close to $1$ (e.g. \cite{Gamburd2002}) so part (i) 
above 
is not vacuous; moreover the example of $\G$ with $\gd$ arbitrarily close
to $1$ given in \cite{Gamburd2002} does in fact contain unipotent elements. 
Part (ii) is not particularly interesting, since one can vastly improve $R$ for a finite index congruence group with classical techniques. We include it here only for comparison.
It is not known whether part (iii)   above is vacuous (even taking $\gd=1$ requires $\gt<5/9\approx.555$ and is outside the reach of Kim-Sarnak), but $R=9$ is the limit of our methods, and coincidentally is precisely the number of factors attained in Brun's original attack \cite{Brun1919} on the twin prime conjecture.

\begin{remark}
Determining membership in $\Or$ amounts to expressing an element in $\G$ as a word in the generators -- not an easy task. We are finding numbers with few prime factors despite having extremely limited knowledge as to which numbers appear!
\end{remark}

\begin{remark}
The choice of the function $f(c,d)=c^2+d^2$ is cosmetic; our methods apply to an arbitrary polynomial $f$, and we plan to detail this generalization in a future publication. That said, our current choice of $f$ is natural, not only historically (indeed the problem of finding primes in sums of two squares dates back to Fermat)  but also because, as we shall  see, it is everywhere unobstructed in the affine linear sieve (akin to looking for primes congruent $1$ modulo $q$ -- there are no obstructions for any $q$).  See Remark \ref{rmkA2}.
\end{remark}

\begin{remark}\label{whybaad}
The requirement in the Main Theorem that $\Gi$ be nontrivial is undesirable. 
Indeed the conclusions should hold without this assumption, and we are currently working to remove it by other methods. 
See Remark \ref{baad}.

More importantly, unipotent elements furnish an affine injection into our orbit, enabling more classical sieve techniques. Precisely, if $\ga,\g\in\G$ and $\ga\neq I$ fixes infinity, then $f\big((0,1)\g\ga^n\big)$ is a quadratic polynomial in $n$. These are known \cite{Iwaniec1978} to contain infinitely many numbers with at most two factors! Moreover by varying $\g$, one can accrue a Zariski dense set of 2-almost primes, in the sense of \cite{BourgainGamburdSarnak2006}. 
 Therefore we state one more immediate corollary which cannot be deduced by ``cheating'' with unipotents:
\end{remark}

\begin{corollary}\label{corcor}
Let $\G,\ \gd$ and $\Or(T,R)$ be as in the Main Theorem. 
Then for $R=1$ we have the following upper bound for the number of primes in $f$:
\benn
|\Or
(T,1)|\ll T^{\gd}/\log T,\quad\quad\text{as $T\to\infty$}.
\eenn
This is off by a constant multiple from the expected asymptotic formula.
\end{corollary}

This paper is organized as follows. In the next section we give background material on  Strong Approximation, the geometry and spectra of infinite co-volume groups, and the weighted linear $\gb$-sieve. In \S 3 we prove a certain Main Identity, which 
shows how to grow the lattice point count at time $T$ from that at small times   via the Laplace operator (this is the key to circumventing an explicit spectral theorem).
 In 
\S 4 we collect  preliminary facts about infinite volume lattice point counts before 
proving the Main Theorem in 
 \S 5. 
Some technical issues are reserved for the Appendices.

\subsection*{Acknowledgements}
I thank my advisors Dorian Goldfeld and  Peter Sarnak for their guidance, encouragement, and inspiration.
 I am grateful to the anonymous referees and Hee Oh for detailed comments and corrections to an earlier draft of this document.


\section{Background Material}

\subsection{Strong Approximation}

Our first ingredients are algebraic in nature. We require the Strong Approximation Theorem of Matthews, Vaserstein, and Weisfeiler \cite{MatthewsVasersteinWeisfeiler1984}, stating that if $L\subset\SL_n(\Z)$ is Zariski dense in $\SL_n(\Z)$ then the projection of $L$ on $\SL_n(\Z_p)$ is dense for all but finitely many primes $p$. 
Recall that the ring of $p$-adic integers $\Z_p$ is the inverse limit of the finite rings $\Z/p^k\Z$.
In particular, this means $L$ is onto $\SL_n(\Z/p\Z)$. Actually 
in dimension two
 this
can be done 
by more elementary methods 
\cite{DavidoffSarnakValette2003}.

We also require Goursat's Lemma which states the following. Let $G_1$, $G_2$ be groups, and let $H$ be a subgroup of $G_1\times G_2$ such that the two projections $p_j: H\rightarrow G_j$,
$j=1,2$
 are surjective. 
Let $N_j$ be the kernel of $p_j$. 
One can identify $N_1$ as a normal subgroup of $G_2$, and $N_2$ as a normal subgroup of $G_1$. Then the image of $H$ in $G_1/N_2\times G_2/N_1$ is the graph of an isomorphism $G_1/N_2\approx G_2/N_1$.
As the group $\PSL_2(\F)$ is simple for a finite field $\F$ satisfying $|\F|>4$, this implies that if $p_1\neq p_2$ are two primes and  $L\subset\SL_2(\Z)$ projects onto both $SL_2(\Z/p_1\Z)$  and $\SL_2(\Z/p_2\Z)$, then $L$ projects onto $SL_2(\Z/p_1 p_2 \Z)$.

We follow 
\cite{BourgainGamburdSarnak2006} in recording these two facts as follows:

\begin{theorem}\label{strongApp}
Let $\G\subset\SL_2(\Z)$ be nonelementary. Then there exists a number $\B$ which is the product of a finite set of ``ramified'' primes  such that if 
 $q=q'q''$ is square-free 
 with $q'|\B$ and $(q'',\B)=1$ then the projection of $\G$ in $\SL_2(\Z/q\Z)$ is 
the product 
 $
G_{q'}\times 
 \SL_2(\Z/q''\Z)$, where $G_{q'}$ is the projection of $\G$ in $\SL_2(\Z/q'\Z)$.
\end{theorem}


\subsection{Infinite-volume geometry and spectral data}\label{delta}

Let $\G\subset\SL_2(\Z)$ be a finitely-generated non-elementary Fuchsian group.
 Then $\G$ acts on the Poincar\'{e} upper half plane $\boldH$ by fractional linear transformations. 
 In dimension two, being finitely generated is equivalent to being geometrically finite, i.e. that the 
 Riemann surface $\Fu=\GbkH$ has finitely many bounding sides \cite{Beardon1983}.

As the action of $\G$ on  $\H$ is discrete, there are no limit points in $\H$. There are however  limit points in the boundary $\widehat{\R}=\R\cup\{\infty\}$. The set of all limit points of $\G$ is called the limit set $\gL=\gL(\G)$. It is a Cantor-like fractal and has some Hausdorff dimension $\gd=\gd(\G)\in[0,1]$. A geometrically finite group has $\gd<1$ if and only if $\Fu$ has infinite hyperbolic volume (i.e. $\G$ has infinite index in $\SL_2(\Z)$).

\begin{figure}
\includegraphics[width=3.2in]{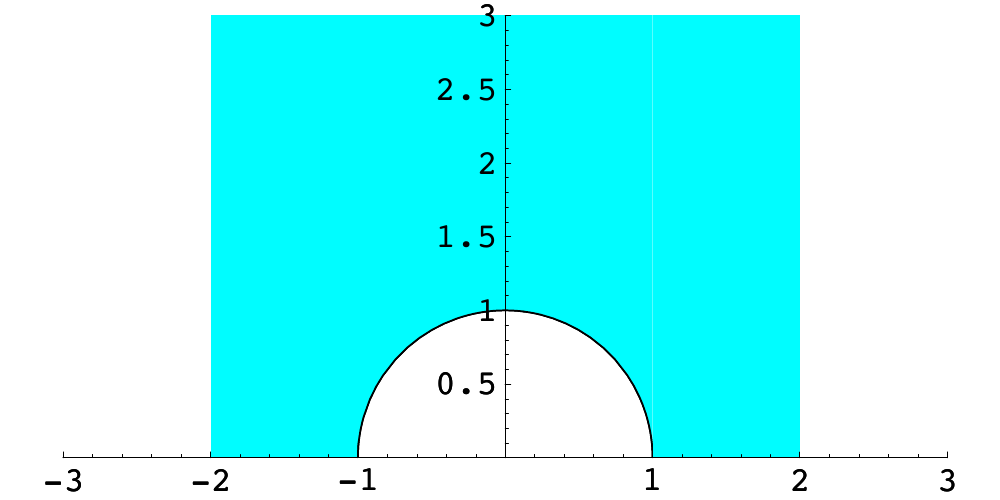}  
\caption{A fundamental domain for $\G=\< T^4, S\>$.}
\label{fig1}
\end{figure}

 If we label the generators of $\SL_2(\Z)$ by $T:z\mapsto z+1$ and \break $S:z\mapsto -1/z$, then a prototypical example of the type of group we have in mind is the Hecke group $\G=\<T^4,S\>$. This group has fundamental domain $\Fu=\{z\in\boldH\mid |z|>1,|\Re(z)|<2\}$, whose vertical strips 
  touch the real line 
at a free boundary
 and
  clearly contribute infinite hyperbolic volume. See Fig. \ref{fig1}. We will assume throughout that 
  $\Fu$ has a cusp at
  infinity, as in this example. 

 As usual, $\Fu$ is equipped with a hyperbolic volume element $dz={dx\,dy\over y^2}$. The positive-definite Laplace-Beltrami operator $\gD=-y^2\left(\ddxt+\ddyt\right)$ of $\Fu$ acts in the space of $C^{\infty}$ functions with compact support $C^{\infty}_K(\Fu)$, and has a unique self-adjoint extension to an unbounded operator on $L^2(\Fu)$. Denote the spectrum of $\gD$ on $L^2(\GbkH)$ by $\Spec(\GbkH)$. 

The spectrum of $\gD$ below $1/4$ consists only of finitely-many point eigenvalues 
and the tempered spectrum contained in  $[1/4, \infty)$ is purely continuous \cite{LaxPhillips1982}. Notice that in finite volume $\gd=1$ and the base eigenvalue is $\gl_0=0$, corresponding to the constant function $\phi_0=1/\sqrt{\vol(\Fu)}$, scaled to have unit $L^2$-norm. Neither does this make sense if $\vol(\Fu)=\infty$ nor is any non-zero constant function  square-integrable in this case.  It follows from the work of Patterson \cite{Patterson1976, Sullivan1984} that  $\gd>1/2$ if and only if 
$$
\Spec(\GbkH)\cap[0,1/4)\neq\O,
$$ 
in which case $\gl_0=\gd(1-\gd)$ is the base eigenvalue of $\gD$. This eigenvalue is isolated, has multiplicity one, and any associated eigenfunction is of constant sign on $\Fu$; in particular we can choose it to be nonnegative.  The r\^{o}le of the constant function (volume) is then played by the base eigenfunction $\phi_0$ which Patterson determined explicitly as  the integral of a Poisson kernel against the so-called Patterson-Sullivan measure $\mu$, supported on the limit set $\gL$:
\be\label{phi0def}
\phi_0(x,y):=\int_{\gL} \left({(t^2+1) y\over (x-t)^2+y^2}\right)^{\gd} d\mu(t).
\ee
If $\Fu$ has a cusp 
(which we assume it does)
 then $\gd>1/2$. 

Let $\G(q)$ denote the principal ``congruence'' subgroup of $\G$ of level $q$,
\be\label{def:Gq}
\G(q):=\{\g\in\G:\g\equiv I(\mod q)\}.
\ee
This is of course still a Fuchsian group of the second kind, but has finite index in $\G$ (and therefore also has the same limit set and Hausdorff dimension -- every Cauchy sequence in $\boldH$ under the action of $\G$ has a corresponding sequence under $\G(q)$ with the same limit point). 
The inclusion $\G(q)\subset\G$ induces the  reverse inclusion 
$$
\Spec(\GbkH)\subset \Spec(\G(q)\bk\boldH)
.
$$ 
In particular this means the base eigenfunction $\phi_0$ is an ``oldform'' on $\G(q)$ (but must be rescaled to have unit $L^2$-norm).

Corresponding to any point eigenvalue $\gl\in\Spec(\G(q)\bk\boldH)$ is an eigenfunction, $\phi\in L^2(\G(q)\bk\boldH)$. 
It may be the case that $\phi$ (which a priori is only $\G(q)$-automorphic) is also automorphic with respect to $\G$.
 In this case we we call $\phi$  an ``oldform'' and  $\gl$ is an ``oldvalue''. 
 In the opposite case we call $\phi$  a ``newform'' and $\gl$ a ``newvalue''. 
Denote by 
$
\Spec(\G(q)\bk\boldH)_{new}
$
the subset of $\Spec(\G(q)\bk\boldH)$ consisting of ``new'' eigenvalues. 

\begin{definition}
We say that $\G$ has a {\bf spectral gap} $\gt\in[1/2,\gd)$ if there exists $\B\in\N$ such that  for  $q=q'q''$ square-free with $q'|\B$ and $(q'',\B)=1$, we have 
\be\label{specGapGt}
\Spec(\G(q)\bk\boldH)_{new}\cap(0,\gt(1-\gt)) \subset \Spec(\G(q')\bk\boldH)_{new}.
\ee
\end{definition}

Note that this definition of the ``spectral gap'' is not
the conventional one (for which see, for example, \cite{EinsiedlerMargulisVenkatesh2007, KelmerSarnak2008}).

Collecting the results in \cite{BourgainGamburd2007,BourgainGamburdSarnak2006} and their extension from prime to square-free of \cite{Gamburd2002} we have:

\begin{theorem}[\cite{Gamburd2002,BourgainGamburd2007,BourgainGamburdSarnak2006}]\label{BGSspec}
Let $\G\subset\SL_2(\Z)$ and
 $\gd>1/2$ be the Hausdorff measure of its limit set.
\begin{enumerate}
\item There exists a spectral gap $\gt\in[1/2,\gd)$ such that
 \eqref{specGapGt} holds with $\B=1$.
\item If $\gd>5/6$ then \eqref{specGapGt} holds with $\gt=5/6$ and 
the number $\B$ 
is precisely the one which appeared in Theorem \ref{strongApp}.
\end{enumerate}
\end{theorem}

The case $\G=\SL_2(\Z)$,  $\B=1$ and $\gt=1/2$  is the celebrated (and unsolved) Selberg $1/4$-Conjecture \cite{Selberg1965}, which in modern parlance is a consequence of  the generalized Ramanujan Conjectures.

We now record the abstract spectral theorem for unbounded self-adjoint 
operators as follows.
\begin{theorem}[Abstract Spectral Theorem]\label{AbsSpecThm}
There exists a spectral measure $\nu$, supported on $\Spec(\G\bk\boldH)$ and a  unitary spectral operator $\hat{\,} :L^2(\G\bk\boldH) \to L^2((0,\infty),d\nu)$ such that:
\begin{enumerate}
\item We have the Abstract Parseval's Identity: for $\phi_1,\phi_2 \in L^2(\G\bk\boldH)$,
\be\label{absPars}
\<\phi_1,\phi_2\>_{L^2(\GbkH,dz)} = \<\hat{\phi_1},\hat{\phi_2}\>_{L^2(\Spec(\GbkH),d\nu)}
\ee
\item The spectral operator \ $\hat{\,}$ \ is diagonal with respect to the Laplacian: for $\phi\in L^2(\GbkH)$ and $\gl\in\Spec(\GbkH)$
\be\label{specDiag}
\hat{\gD \phi}(\gl)=\gl\,\hat{\phi}(\gl).
\ee
\item If $\gl$ is a point eigenvalue of multiplicity one with associated $L^2$ eigenfunction $\phi_{\gl}$ of unit norm, then for any $\phi\in L^2(\GbkH)$,
\be\label{phigl}
\hat{\phi}(\gl) = \<\phi,\phi_{\gl}\>.
\ee
\end{enumerate}
\end{theorem}

\subsection{Combinatorial Sieve}\label{combsieve}

Let $\A$ be a 
sequence of non-negative real numbers $\{a_n\}_n$ of which all but finitely many 
are zero.
For $R\ge1$ let  $Z_R$ denote the set of positive integers with at most $R$ prime divisors.
The main objective in sieve theory is to determine lower bound  estimates for
 \be\label{sumAnS}
\sum_{n\in Z_R}a_n
\ee
 given knowledge of how $\A$ is distributed along each of the arithmetic progressions $0(\mod q)$ for square-free $q$. 
In the following setup of the sieve, there are many parameters. Their heuristic  meaning is as follows. 
 
 Let $q\ge1$ be square-free and collect the elements of $\A$ whose index is divisible by $q$ via 
\be\label{Aqs}
\A_q:=\{a_n\in\A\mid n\equiv0(q)\}.
\ee
Of course $\A_1=\A$. 
%
The parameter   $T$ 
is the cut-off point after which all $a_n$'s are zero.
 The parameter $\B\in\N$ is the product of a fixed finite set of ``ramified'' primes. We will decompose $q$ into the $\B$-part and the rest: $q=q'q''$ with $q'|\B$ and $(q'',\B)=1$. The parameter $\X$ is an approximation to $|\A|:=\sum_n a_n$, and $\gw:\N\to[0,1]$ represents the local density at $q$. Thus $\gw(q)\X$ is an approximation to 
 $$
 |\A_q|=\sum_{n\equiv0(q)}a_n.
 $$ 
 At the $\B$-part, there is a lower order term, $\X_{q'}$ such that 
 $$
 r(q):=|\A_q|-\gw(q)(\X+\X_{q'})
 $$ 
 is an error term. The error is  small on average: there is a sieving level $Q$ such that the total contribution from the error terms up to $Q$ is  a power savings off the main term $\X$.

Precisely, we require the following conditions: 
\renewcommand{\labelenumi}{(S\theenumi)}
\begin{enumerate}
\item  $\X\asymp |\A|:=\sum_n a_n$
and
$\B\in\N$ is a fixed natural number. There is a parameter  $T>1$ such that
\be\label{aNwithTbound}
a_n = 0 \text{ if } n>T.
\ee
\item \label{gwsat} The function $\gw:N\to[0,1]$ satisfies $\gw(1)=1$ and $\gw(q)<1$ for $q>1$. 
Moreover, $\gw$ is multiplicative away from $\B$. 
By this we mean that
for any $q$  square-free, write $q=q'q''$ with $q'|\B$ and $(q'',\B)=1$. Then $\gw(q)=\gw(q')\prod_{p|q''}\gw(p)$.
\item\label{localDensBound}  There exists a fixed constant $K<\infty$ such that for $2\le v\le z$, we have  the following 
local density 
bound: 

$\prod_{{v\le p \le z}
\atop{p\nmid\B}
} (1-\gw(p))^{-1} \le  \left(\frac{\log z}{\log v}\right)\left(1+{K\over\log v}\right).
$
\item \label{extraTerms}
For any divisor $q'|\B$ let $\X_{q'}$ satisfy $\X_{q'}\ll\X^{1-\eta}$
for some
$\eta>0$.
\item  \label{rqErr}
 Define 
$r(q):=|\A_q|-\gw(q)\left(\X + \X_{q'}\right)$
and assume that for some parameter $Q>1$ and $\gep>0$,
$$
\sum_{q\le Q} |r(q)| \ll \X^{1-\gep}.
$$
\end{enumerate}

The following Theorem is a simple consequence of Theorem 3.7 on page 63 in \cite{Iwaniec1996}. We derive it from the original in Appendix \ref{AxSieve}.

\begin{theorem}\label{sieve}
Let $\A$ be as described above. Then 
\be\label{eq:Sz}
 \sum_{n\in Z_R} a_n \asymp \X \prod_{{p\le Q
}\atop{p\nmid\B} }(1-\gw(p)),
\ee
for any $R$ satisfying 
\be\label{Rbound}
R>2\log T/\log Q.
\ee

\end{theorem}

In our application, Goursat's Lemma demonstrates the multiplicativity of the local density $\gw$ in (S\ref{gwsat}), and  Strong Approximation allows us to compute $\gw$  explicitly and verify (S\ref{localDensBound}). 
The finite collection of primes comprising $\B$ 
which may have exceptional eigenvalues in Theorem \ref{BGSspec}
are directly responsible for the extra terms $\X_{q'}$ in (S\ref{extraTerms}) but their contribution is harmless to the inclusion-exclusion of the sieve.


\section{The Main Identity}\label{mainId}
Let $\G\subset\SL_2(\Z)$.
Recall the classical fact that for $\g=\mattwo{*}{*}{c}{d}\in\G$ we have
\benn
\Im(\g z)={\Im z\over |cz+d|^2}.
\eenn
Fix $T>1$ and define the characteristic function
\be\label{chiTdef}
\chi_T(z):=\twocase{}{1}{if $\Im(z)>1/T$}{0}{otherwise,}
\ee
so that 
\benn
\chi_T(\g\i)=\twocase{}{1}{if $c^2+d^2<T$}{0}{otherwise.}
\eenn
Let $\Gi\subset\G$ be the maximal 
subgroup stabilizing infinity; clearly $\chi_T$ is $\Gi$-invariant.

We average $\chi_T$ over the group $\G$ 
\be\label{FTdef}
F_T(z):=\sum_{\g\in\Gi\bk\G}\chi_T(\g z),
\ee
so that we can recover the lattice point count \eqref{def:OT} via
\benn
F_T(\i)=|\Or(T)|.
\eenn
Clearly $F_T$ is $\G$-invariant, so is well defined as a function on $\GbkH$. 
\begin{lemma}
$F_T\in L^2(\GbkH)$ if and only if $\Gi$ is nontrivial.
\end{lemma}
\begin{proof}
If $\Gi$ is nontrivial,
let $N$ be the 
finite 
interval corresponding to the width of the cusp at infinity. 
Otherwise set $N=\R$. 
 By unfolding we have
\beann
\|F_T\|^2&=&\int_{\GbkH} F_T(z)F_T(z)dz = \int_{\GbkH} \sum_{\g\in\Gi\bk\G}\chi_T(\g z) F_T(z)dz\\
&=& \int_{\GibkH}\chi_T(z)F_T(z) dz= \int_{1/T}^{\infty}\int_N F_T(z) dx{ dy\over y^2}.
\eeann
For $\Im(z)$ sufficiently large, only the term $\g=I$ contributes to the sum \eqref{FTdef}, and so $F_T(z)=1$. It is now clear that the above integral converges if and only if $N\neq\R$.
\end{proof}

Instead of accessing $F_T$ directly, we will seek an identity which ``grows'' the count at time $T$ from the count at time $T=1$ (and at some other time, $T=b$). Let $\phi\in L^2(\GbkH)$ be an eigenfunction of the hyperbolic Laplace-Beltrami operator $\gD=-y^2\left(\ddxt+\ddyt \right)$ with eigenvalue $\gl=s(1-s)<1/4$. By unfolding the following inner product and using the fact that the constant (or any) Fourier coefficient of $\phi$ satisfies the same differential equation as $\phi$ itself, we have
\beann
\<F_T,\phi\>&=& \int_{\GbkH} F_T(z)\bar{\phi}(z)dz
= \int_{\GibkH} \chi_T(z)\bar{\phi}(z)dz \\
&=& \int_{1/T}^{\infty} \left(\int_N \bar{\phi}(z)dx\right){dy\over y^2}
=\int_{1/T}^{\infty} (\ga y^s+\gb y^{1-s}){dy\over y^2} \\
&=& A_{\phi} T^s + B_{\phi}T^{1-s},
\eeann
where $\ga,\ \gb,\ A_{\phi},$ and $B_{\phi}$ are some constants depending on the eigenfunction $\phi$.
(As $F_T$ is like a truncated Eisenstein series, this calculation is just the incomplete Mellin transform of the constant term of $\phi$.) 

Following the methodology of Selberg \cite{Selberg1956}, we seek an identity which depends only on the eigenvalue $\gl$ or $s$ but not on the eigenfunction $\phi$. Reformulate the above identity as a dot product of vectors
\be\label{ABT}
\<F_T,\phi\>=A_{\phi}T^s+B_{\phi}T^{1-s}  = (T^s, T^{1-s})\vectwo{A_{\phi}}{B_{\phi}}
\ee
and write it at time $T=1$ and some other time $T=b>1$:
\beann
\<F_1,\phi\>&=&A_{\phi}+B_{\phi},\\
\<F_b,\phi\>&=&A_{\phi}b^s+B_{\phi}b^{1-s},
\eeann
or
\beann
\vectwo{\<F_1,\phi\>}{\<F_b,\phi\>} & =& \mattwo{1}{1}{b^s}{b^{1-s}}\vectwo{A_{\phi}}{B_{\phi}}.
\eeann
Notice that the matrix on the right hand side only depends on the eigenvalue $s$. We multiply both sides by the inverse 
matrix
\benn
\vectwo{A_{\phi}}{B_{\phi}}= {1\over b^{1-s}-b^s}\mattwo{b^{1-s}}{-1}{-b^s}{1} \vectwo{\<F_1,\phi\>}{\<F_b,\phi\>},
\eenn
and insert this into \eqref{ABT}:
\bea \nonumber
\<F_T,\phi\>
& =& (T^s, T^{1-s}) {1\over b^{1-s}-b^s}\mattwo{b^{1-s}}{-1}{-b^s}{1} \vectwo{\<F_1,\phi\>}{\<F_b,\phi\>}\\ \label{FKL}
&=&K_T(s) \<F_1,\phi\> + L_T(s)\<F_b,\phi\>,
\eea
where
\bea\label{KLdef}
K_T(s)&=&{T^s b^{1-s}-T^{1-s}b^s\over b^{1-s}-b^s},\\ \nonumber
L_T(s)&=&{T^{1-s}-T^s\over b^{1-s}-b^s}
\eea
are functions which only depend on the eigenvalue, $\gl=s(1-s)$.

Notice that for $T>1$, $s\in(1/2,1]$,  and  $b$ fixed
\be\label{KLTbound}
K_T(s),L_T(s)\ll T^s.
\ee
For $s=\foh+\i t$ we have
\beann
K_T(s)&=&T^{1/2}{\sin(t\log T/\log b)\over \sin(t\log b)},\\
L_T(s)&=&\left({T\over b}\right)^{1/2}{\sin(t\log T)\over \sin(t\log b)}.\\
\eeann
Choosing $b=b(T)$ such that ${\log T\over\log b}\in\Z$ ensures that the functions above are 
\[
\ll T^{1/2} \log T. 
\]
For example, we can take
\be\label{bdef}
b=\exp\left({\log T\over \lceil \log T\rceil}\right)\in(1,e],
\ee
where $\lceil \cdot \rceil$ is the ceiling function, returning the smallest integer not less than its argument. In particular, 
\[
b<3.
\]

As $\gl=s(1-s)$, we abuse notation by writing 
\[
K_T(\gl) \text{ and } L_T(\gl)
\] 
in place of 
\[
K_T(s) \text{ and }L_T(s).
\] 
This should cause no confusion.
Just as one can exponentiate a matrix, one can define a function of a differential operator, which is itself a differential operator. 
So as $\phi$ is an eigenfunction satisfying $\gD\phi=\gl\phi$, 
we have
\benn
K_T(\gD)\phi = K_T(\gl) \phi,
\eenn
where the left hand side acts by differentiation  and the right hand side is multiplication by a function of the eigenvalue. The same holds for $L_T$. 

Since the Laplacian $\gD$ is self-adjoint, we have from \eqref{FKL}
\beann
\<F_T,\phi\> &=&  K_T(\gl)\<F_1,\phi\>+L_T(\gl)\<F_b,\phi\>\\
&=&\<F_1,K_T(\gl)\phi\>+\<F_b,L_T(\gl)\phi\>\\
&=&\<F_1,K_T(\gD)\phi\>+\<F_b,L_T(\gD)\phi\>\\
&=&\left<K_T(\gD)F_1,\phi\right> + \left<L_T(\gD)F_b,\phi\right>\\
&=&\left<K_T(\gD)F_1 + L_T(\gD)F_b,\phi\right>.
\eeann
This identity holds for any eigenfunction $\phi$ in the point spectrum, and so it should be the case that
\be\label{mmainid}
F_T(z) = K_T(\gD) F_1(z) + L_T(\gD) F_b(z)
\ee
holds in general. Notice that the $T$ dependence has been entirely removed from $F$ on the right hand side and only appears in the spectral operators $K_T$ and $L_T$! It should now be clear  what we mean by gathering information about $F$ at time $T$ using only the Laplacian and $F$ at small times. We reserve the rest of this section for the proof of this Main Identity.

\begin{theorem}[Main Identity]\label{structThm}
For fixed $T\ge1$ let $F_T$ be defined by \eqref{FTdef}. 
Then there exists a number $b<3$ and functions $K_T$ and $L_T$ satisfying:
\benn
L_T(\gl)\ll K_T(\gl)\ll T^s
\eenn
for $\gl=s(1-s)$ and $s\in(1/2,1]$, and 
\benn
L_T(\gl)\ll K_T(\gl)\ll T^{1/2}\log T
\eenn
for $\gl\ge1/4$, 
such that  \eqref{mmainid} holds for almost every $z$.
\end{theorem}

The argument above \eqref{mmainid} proves this identity along the point spectrum. If we had an explicit spectral theorem in infinite volume, we would just need to carry out similar  computations to prove this identity on  the continuous spectrum (for a finite co-volume group, this is easily achieved via the Eisenstein series). Instead, we will use 
ideas from ``almost'' eigenfunctions and perturbation theory  for the proof, which will occupy the remainder of this section.

We begin with a technical lemma. Note that from now on, $\phi$ is not assumed to be an eigenfunction of $\gD$, only that $\phi\in L^2(\GbkH)$. 

\begin{lemma}\label{lem:ODE}
For any $\phi\in L^2(\GbkH)$ and any $\gl=s(1-s)\ge0$, $\gl\neq1/4$ there exist constants $A$ and $B$ such that
\be\label{ODE}
\<F_T, \phi\>=A\, T^s+B\, T^{1-s}+O_{\gl, T,\G} \left(\|(\gD-\gl)\phi\|\right).
\ee
If $\gl=1/4$ then
\benn
\<F_T, \phi\>=A\, T^{1/2}+B\, T^{1/2}\log T+O_{\gl, T,\G} \left(\|(\gD-\gl)\phi\|\right).
\eenn
\end{lemma}

This is an explicit computation, applying the method of variation of parameters and estimating the inhomogeneous component
(which relies on $\Gi$ being non-trivial). 
We present the details in Appendix \ref{technic}.

\begin{remark}
One should think of $\phi$ as an ``almost'' eigenfunction with eigenvalue $\gl$. 
Then the error term in \eqref{ODE} should be small. This is a heuristic only; the argument applies for general $\phi$ and $\gl$, and we will soon see that such an application is necessary.
\end{remark}

The following proposition says that the difference of the right and left hand sides of \eqref{ODE} has no correlation with any almost eigenfunction.

\begin{proposition}\label{prop:G}
Let $T,\ b,$ and $F_T$ be as above and set 
\be\label{GTdef}
G_T=F_T-K_T(\gD)F_1-L_T(\gD)F_b.
\ee 
Then for any  $\phi\in L^2(\GbkH)$ and any $\gl\ge0$ we have
\be\label{eq:Gineq}
\<G_T,\phi\> \ll_{\gl,T}\ \| (\gD-\gl)\phi\|. 
\ee
\end{proposition}

\begin{proof}

Fix an arbitrary $\phi\in L^2(\GbkH)$ and $\gl=s(1-s)\ge0,$ and let $G_T$ be defined by \eqref{GTdef}.
Assume $\gl\neq1/4$ (the computation in the the case of $\gl=\frac14$ is similar).
Consider the following trivial identity obtained by adding and subtracting identical terms:
\bea
\nonumber \<G_T,\phi\> & = & \< F_T, \phi\> - (A_{\phi} T^s + B_{\phi} T^{1-s}) \\
\label{RowTwo}  & & - \left( \< K_T(\gD)F_1,\phi\> - K_T(\gl)\< F_1,\phi\>  \right) \\
\nonumber  & & - \left( \< L_T(\gD)F_b,\phi\> - L_T(\gl)\< F_b,\phi\> \right) \\
\nonumber  & & - \left( K_T(\gl)\< F_1,\phi\> - K_T(\gl)(A_{\phi}+B_{\phi})\right) \\
\nonumber  & & - \left(L_T(\gl)\< F_b,\phi\> - L_T(\gl) (A_{\phi}b^s+B_{\phi}b^{1-s})\right) \\
\nonumber  & & + (A_{\phi}T^s + B_{\phi}T^{1-s}) - K_T(\gl) (A_{\phi} + B_{\phi}) - L_T(\gl) (A_{\phi} b^s + B_{\phi} b^{1-s}).
\eea
The bottom row is zero by construction of $K_T$ and $L_T$ in \eqref{KLdef}. The top row is $\ll_{\gl,T} \|(\gD-\gl)\phi\|$ by Lemma \ref{lem:ODE}, as are the fourth and fifth rows. It remains to understand the second and third rows. 

Let $\hat{\phi}$ denote the spectral transform of $\phi$ and let $\hat{F_T}$ be the spectral transform of $F_T$.  As $\gl$ is fixed, let $\gl'$ be in $\Spec(\GbkH)$.
By the Mean Value Theorem and Cauchy-Schwarz, the second row is
\beann
\eqref{RowTwo}&=&
\< K_T(\gD)F_1, \phi\> - K_T(\gl)\<F_1, \phi\> \\
& = & \< K_T\cdot \hat{F_1},\widehat{\phi}\> - K_T(\gl)\< \hat{F_1},\widehat{\phi}\> \\
 & = & \int_{\Spec(\GbkH)} (K_T(\gl')-K_T(\gl))  {\hat{F_1}}(\gl')\widehat{\bar{\phi}}(\gl')d\nu(\gl')\\
 & \ll_{\gl,T} & \int_{\Spec(\GbkH)} (\gl'-\gl) \widehat{\bar{\phi}}(\gl') {\hat{F_1}}(\gl')  d\nu(\gl')\\ 
  & \ll & \|(\gD-\gl)\phi\| \|F_1\|.
\eeann

The calculation for the third row is identical and we are done.

\end{proof}

Finally, we show that a function which is uncorrelated to any almost eigenfunction is zero.

\begin{proof}[Proof of the Main Identity.]

We aim to show that  $G_T$
defined in \eqref{GTdef}
 vanishes almost everywhere. 

Since \eqref{eq:Gineq} holds for any $\phi$, we  are free to choose our $\phi$. Equivalently we may choose its spectral transform $\widehat{\phi}$, so we make the following construction. Fix $\vep>0$ and fix an arbitrary $\gl>0$. Let
\be\label{phiHat}
\widehat{\phi}(\gl') := \twocase{}{\widehat{G}_T(\gl')}{if $\gl'\in(\gl-\vep,\gl+\vep)$}{0}{otherwise,}
\ee
where $\widehat{G}_T$ is the spectral transform of $G_T$.
Inserting \eqref{phiHat} into  Abstract Parseval's Theorem \eqref{absPars} we have:
\bea
\nonumber \<G_T,\phi\> &= & \< \widehat{G}_T, \widehat{\phi}\> = \int_{\Spec(\GbkH)}\widehat{G}_T(\gl')\widehat{\bar{\phi}}(\gl')d\nu(\gl')\\
\label{GTRHS} & = &  \int_{\gl-\vep}^{\gl+\vep} |\widehat{G}_T(\gl')|^2 d\nu(\gl'),
\eea
and
\bea 
\nonumber \|(\gD-\gl)\phi\| & = & \left( \int_{\Spec(\GbkH)} \left|(\gl'-\gl)\widehat{\phi}(\gl')\right|^2 d\nu(\gl')\right)^{1/2} \\
\nonumber  & = & \left( \int_{\gl-\vep}^{\gl+\vep}  | \gl' - \gl |^2 \left| \widehat{G}_T(\gl')\right|^2 d\nu(\gl') \right)^{1/2} \\
\label{GTLHS} & \le & \vep \left( \int_{\gl-\vep}^{\gl+\vep} \left| \widehat{G}_T(\gl')\right|^2 d\nu(\gl') \right)^{1/2}. 
\eea

Inserting \eqref{GTRHS} and \eqref{GTLHS} into  \eqref{eq:Gineq}, we have:
\benn
\int_{\gl-\vep}^{\gl+\vep} |\widehat{G}_T(\gl')|^2 d\nu(\gl') \ll_{\gl,T} \vep \left( \int_{\gl-\vep}^{\gl+\vep} |\widehat{G}_T(\gl')|^2 d\nu(\gl') \right)^{1/2}.
\eenn
If the left side is zero, we are done. If not, we have for an arbitrary $\gl$:
\benn
\int_{\gl-\vep}^{\gl+\vep} |\widehat{G}_T(\gl')|^2 d\nu(\gl') \ll_{\gl,T} \vep^2.
\eenn
Let 
$f(\gl)
=\int_{\gl' < \gl} |\widehat{G}_T(\gl')|^2d\nu(\gl')$.
 Then $f$ is everywhere continuously differentiable with
\benn
f'(\gl)=\lim_{\vep\to0}{f(\gl+\vep)-f(\gl-\vep)\over 2\vep}= \lim_{\vep\to0} {1\over2\vep}\int_{\gl-\vep}^{\gl+\vep} |\widehat{G}_T(\gl')|^2 d\nu(\gl')=0.
\eenn
So as $f'=0$ uniformly and $f(0)=0$, we have shown that $f\equiv0$, and  $G_T=0$ a.e.

This concludes the proof of the Main Identity.
\end{proof}


\section{Preliminaries}

\subsection{Sums over $w_T$}

Once again let $\G\subset\SL_2(\Z)$ and let $\Gi\subset\G$ be the stabilizer of infinity. Fix $\vep>0$ and let 
\be\label{psidef}
\psi=\psi_{\vep}\in L^2(\GibkH)
\ee 
be an $\vep$-approximation to the identity about $z_0=\i$. By this we ask that $\psi$ 
\renewcommand{\labelenumi}{(\theenumi)}
\begin{enumerate}
\item be smooth, nonnegative, 
\item have total mass one $\int_{\GibkH}\psi=1$, and 
\item be supported in a small neighborhood about $z_0=\i$ (a ball of radius $\vep/10$ will suffice).
\end{enumerate}
Recall the characteristic function $\chi_T$ from \eqref{chiTdef}. Define the function $w_T=w_{T,\vep}:\Gi\bk\G\to[0,1]$ by
\be\label{wTdef}
w_T(\g):=\int_{\GibkH} \chi_T(z) \psi_{\vep}(\g z) dz.
\ee
For $\g\in\G$ having bottom row $(c,d)$ it is easy to see that 
\be\label{dubyabound}
w_{T
,\vep
}(\g)=\twocase{}{1}{if $c^2+d^2<T/(1+\vep)$}{0}{if $c^2+d^2>T/(1-\vep)$.}
\ee

We first prove   the following lemma. Let  $\Xi$ be a subgroup of $\G$ having finite index in $\G$. Then both groups  $\Xi$ and $\G$ have the same limit set, Hausdorff dimension $\gd$ and base eigenvalue $\gl_0=\gd(1-\gd)$. Assume further that $\Xi$ contains $\Gi$.

\begin{lemma}\label{GG0Game}
Let $\Gi\subset \Xi\subset \G$ with $[\G:\Xi]<\infty$ and $w_T$ be defined by \eqref{wTdef}.  
Let $\gl_1=s_1(1-s_1)$ be the first eigenvalue above the base in $\Spec(G\bk\boldH)$, and assume $\gl_1<1/4$.
 Then for any fixed $\g\in\G$, there exists a constant $c_{\G,\vep}>0$ depending only on $\G$ and $\vep$ such that
\be\label{gsum}
\sum_{\xi\in\Gi\bk \Xi}w_T(\xi\g)={c_{\G,\vep} T^{\gd}\over [\G:\Xi] }+ O\left({1\over\vep}
T^{s_1}\right),
\ee
where 
  the implied constant depends on $\G$ but not $\Xi$, $\g$ or $\vep$. Moreover, 
  \be\label{cGamma}
  c_{\G,\vep}=c_\G (1+O(\vep)),
  \ee
  where $c_\G>0$ does not depend on $\vep$.

If $\gl_1\ge1/4$, replace the error term in \eqref{gsum} by $O({1\over\vep}T^{1/2}\log T)$.
\end{lemma}
\begin{proof}

Throughout we will suppress dependence on $\vep$ until it is convenient. 

Let $F_T^\Xi$ be the function defined on $\Xi\bk\boldH$ which is the average of the characteristic function $\chi_T$ over the group $\Xi$:
\be\label{FTGdef}
F_T^\Xi(z):=\sum_{\xi\in\Gi\bk \Xi}\chi_T(\xi z).
\ee
Similarly average $\psi$ over 
$\Xi$ to get
\be\label{Psidef}
\Psi^\Xi_{\vep,z_0}(z) :=\sum_{\xi\in\Gi\bk \Xi}\psi_{\vep,z_0}(\xi z),
\ee
which is  an $\vep$-approximation to the identity about $z_0=\i$ in $L^2(\Xi\bk\boldH)$. 
Notice that $\psi_{z_0}(\g^{-1}z)=\psi_{\g z_0}(z)$, where the latter function   is an approximation to the identity about $\g z_0$.

 Replace $\g$ in \eqref{gsum} by $\g^{-1}$ for convenience.
Input the definition \eqref{wTdef} into the left hand side of \eqref{gsum} and repeatedly unfold and refold the integrals:
\bea \nonumber 
\sum_{\xi\in\Gi\bk \Xi}w_T(\xi\g^{-1}) 
&=&\sum_{\xi\in\Gi\bk \Xi} \left(\int_{\GibkH} \chi_T(z) \psi_{z_0}(\xi \g^{-1}z) dz\right)\\ \nonumber
&=&\sum_{\xi\in\Gi\bk \Xi} \left(\int_{\GibkH} \chi_T(z) \psi_{\g z_0}(\xi z) dz\right)\\ \nonumber
&=&\int_{\GibkH} \chi_T(z)\left(\sum_{\xi\in\Gi\bk \Xi}\psi_{\g z_0}(\xi z)\right)  dz\\ \nonumber
&=&\int_{\GibkH} \chi_T(z) \Psi^\Xi_{\g z_0}(z) dz\\ \nonumber
&=&\sum_{\xi\in\Gi\bk \Xi}\left( \int_{\Xi\bk\boldH} \chi_T(\xi z) \Psi^\Xi_{\g z_0}(z) dz\right)\\ \nonumber
 &=& \int_{\Xi\bk\boldH} F_T^\Xi(z) \Psi^\Xi_{\g z_0}(z) dz\\ \label{eq:innerprodH}
 &=&\<F^\Xi_T,\Psi^\Xi_{\vep,\g z_0}\>_\Xi,
\eea
where
 the inner product above $\<\cdot,\cdot\>_\Xi$ is with respect to the Hilbert space $L^2(\Xi\bk\boldH)$. The above exchanges of summation and integration are justified since everything in sight is  nonnegative and convergent -- the sum on the left hand side of \eqref{eq:innerprodH} has finitely many terms by \eqref{dubyabound}.

By Abstract Parseval's Theorem \eqref{absPars},
\bea\label{Ginnerprod}
\<F^\Xi_T,\Psi^\Xi_{
\g z_0}\>_\Xi &=& \<\hat{F}^\Xi_T,\hat{\Psi}^\Xi_{
\g z_0}\>_{\Spec(\Xi\bk\boldH)} \\ \nonumber
&=& \hat{F}^\Xi_T(\gl_0)\hat{\Psi}^\Xi_{
\g z_0}(\gl_0) + \int_{\Spec(\Xi\bk\boldH)-\{\gl_0\}} \hat{F}^\Xi_T(\gl)\hat{\Psi}^\Xi_{
\g z_0}(\gl) d\nu_\Xi(\gl),
\eea
where $\nu_\Xi$ is the spectral measure on $\Spec(\Xi\bk\boldH)$ (see Theorem \ref{AbsSpecThm}).
By \eqref{phigl} and multiplicity one of the base eigenvalue $\gl_0=\gd(1-\gd)$ we have
\bea\label{FThatP}
\hat{F}^\Xi_T(\gl_0) &=& \<F^\Xi_T,\phi_0^\Xi\>_\Xi,\\ \label{PsiHat}
\hat{\Psi}^\Xi_{
\g z_0}(\gl_0) &=& \<\Psi^\Xi_{
\g z_0},\phi_0^\Xi\>_\Xi,
\eea
where $\phi_0^\Xi$ is the Patterson-Sullivan base eigenfunction in $L^2(\Xi\bk\boldH)$, normalized to have unit norm. Recall this function is real and nonnegative. It is elementary to verify that 
\be\label{phi0G}
\phi_0^\Xi={1\over \sqrt{[\G:\Xi]}} \phi_0^{\G}. 
\ee
Unfolding  \eqref{FThatP} and inserting \eqref{phi0G}, we have
\bea\nonumber
\hat{F}^\Xi_T(\gl_0) &=&  \<F^\Xi_T,\phi_0^\Xi\>_\Xi\\ \nonumber
&=&  \int_{\Xi\bk\boldH} F^\Xi_T(z)  {\phi}_0^\Xi(z)dz\\ \nonumber
&=&  \int_{\Gi\bk\boldH} \chi_T( z)   {\phi}_0^\Xi(z)dz\\ \label{FGTeq}
&=& {1\over \sqrt{[\G:\Xi]}} \int_{1/T}^{\infty} \left(\int_N   {\phi}_0^{\G}(z)dx\right){dy\over y^2},
\eea
where $N$ is a finite interval corresponding to the width of the parabolic at infinity in $\G\bk\boldH$. The key step here is that we have assumed $\Gi\subset \Xi\subset \G$, so $\Xi$ contains the same maximal unipotent subgroup as $\G$. Since $\Gi$ is assumed to be nontrivial, the Hausdorff dimension of the limit set of $\G$ satisfies 
\benn
\gd>1/2. 
\eenn
Therefore $\phi_0$ is an eigenfunction of the Laplacian $\gD$ with eigenvalue $\gl_0=\gd(1-\gd)$, the constant term in its Fourier expansion inherits the same differential equation. Therefore the inner integral in \eqref{FGTeq} is  $\ga y^{\gd} + \gb y^{1-\gd}$, where $\ga$ and $\gb$ are 
some periods (constants depending on $\phi_0$). 
Inserting this into \eqref{FGTeq} and computing the elementary integral, we have
\bea
\nonumber
\hat{F}^\Xi_T(\gl_0) &=&{{1}\over \sqrt{[\G: \Xi]}}\left(  {\gb \over \gd} T^{\gd} +{\ga \over 1-\gd} T^{1-\gd} \right)\\
\label{FTGequals}
&=&  {c'\over \sqrt{[\G: \Xi]}} T^{\gd} + O(T^{1/2}),
\eea
where $c'=\gb/\gd>0$ since $\phi_0$ is nonnegative. Note that $c'$ only depends on $\G$ and not on $\Xi$ or $\vep$ (which has been suppressed until now).

Returning to \eqref{PsiHat},  use \eqref{phi0G} to define $c_{\G,\vep}$ by
\be\label{cGis}
{c_{\G,\vep}
\over {[\G:\Xi]}}
={c'\over \sqrt{[\G: \Xi]}}\<\Psi^\Xi_{\vep,\g z_0},\phi_0^\Xi\>_\Xi
={c'\over {[\G: \Xi]}}\<\Psi^\Xi_{\vep,\g z_0},\phi_0^\G\>_\Xi
,
\ee
as $\phi_0$ is $\Xi$-invariant. Using  \eqref{cGis}, we see that the main term in \eqref{Ginnerprod} coincides with the corresponding term in \eqref{gsum}. 

For more precise information on $c_{\G,\vep}$, in particular its dependence on $\vep$, use the Mean Value Theorem and \eqref{phi0G} to get
\bea
\<\Psi^\Xi_{\vep,\g z_0},\phi_0^\G\>_\Xi
 &=& \int_{\Xi\bk\boldH} \Psi^\Xi_{\vep,\g z_0}(z)\phi_0^\G(\g z_0)dz \\ \nonumber
 & &+ \int_{\Xi\bk\boldH} \Psi^\Xi_{\vep,\g z_0}(z)\left(\phi_0^\G(z)- \phi_0^\G(\g z_0)\right)dz \\ \nonumber
 &=& \phi_0^\G(\g z_0)
 + O\left(
   \sup_{w\in B_{\vep}(\g z_0)}\phi_0'(w) \cdot \vep \right) \\ 
\label{PsiGeq}
 &=& 
 \phi_0^{\G}( z_0) 
 + O\left( \vep\right),
\eea
since $\phi_0^{\G}$ is $\G$-automorphic and $\int \Psi^\Xi =1$. Thus 
  \eqref{PsiGeq}, together with \eqref{cGis}, verifies \eqref{cGamma}, where explicitly, 
  $$
  c_\G= c' \phi_0^\G(\i).
  $$

Returning to the rest of the spectrum in \eqref{Ginnerprod}, we apply the Main Identity \eqref{mmainid} (valid for arbitrary $\G$, in particular for $\G=\Xi$) to the error term:
\bea
\label{Errrr}Err&=& \int_{\Spec(\Xi\bk\boldH)-{\gl_0}} \hat{F}^\Xi_T(\gl)\hat{\Psi}^\Xi_{\vep,\g z_0}(\gl) d\nu_\Xi(\gl)\\
\nonumber&=& \int_{\Spec(\Xi\bk\boldH)-{\gl_0}} \hat{K_T(\gD) F_1^\Xi}(\gl)\hat{\Psi}^\Xi_{\vep,\g z_0}(\gl) d\nu_\Xi(\gl)\\
\nonumber&&+\int_{\Spec(\Xi\bk\boldH)-{\gl_0}} \hat{ L_T(\gD) F_b^\Xi}(\gl)\hat{\Psi}^\Xi_{\vep,\g z_0}(\gl) d\nu_\Xi(\gl),
\eea
where $b<3$.
Assume $\gl_1<1/4$.
By \eqref{specDiag}, $\hat{K_T(\gD) F_1^\Xi}(\gl)=K_T(\gl) \hat{F_1^\Xi}(\gl)$, so together with Cauchy-Schwarz 
and the bound \eqref{KLTbound} we have
\be\label{almostThere}
\int_{\Spec(\Xi\bk\boldH)-{\gl_0}} \hat{K_T(\gD) F_1^\Xi}(\gl)\hat{\Psi}^\Xi_{\vep,\g z_0}(\gl) d\nu_\Xi(\gl) \ll
T^{s_1} \|F_1^\Xi\|_\Xi \|\Psi^\Xi_{\vep}\|_\Xi,
\ee
where $\|\cdot\|_\Xi$ is the norm on $L^2(\Xi\bk\boldH)$.
Clearly the inclusion $\Xi\subset\G$ and positivity of $\chi_T$ implies the pointwise bound
\benn
F_T^\Xi(z)=\sum_{\xi\in\Gi\bk \Xi}\chi_T(\xi \,z) \le \sum_{\g\in\Gi\bk\G}\chi_T(\g z)=F_T^{\G}(z).
\eenn
Applying the pointwise bound  directly 
gives
\benn
\|F_T^\Xi\|_\Xi\le \|F_T^{\G}\|_\Xi = \sqrt{[\G:\Xi]}\|F_T^{\G}\|_{\G},
\eenn
but this is not good enough for us -- we will lose information in the sieve! Instead, we can exploit the positivity of $F_T^\Xi$ to unfold the $L^2$ norm with respect to one copy of $F_T^\Xi$, apply the pointwise bound to the other copy, and refold again:
\bea
\nonumber\|F_T^\Xi\|_\Xi^2 &=& \int_{\Xi\bk\boldH} \left|F_T^\Xi(z)\right|^2dz = \int_{\Xi\bk\boldH} \left(\sum_{\xi\in\Gi\bk \Xi}\chi_T(\xi\, z)\right)F_T^\Xi(z)dz \\
\nonumber&=& \int_{\Gi\bk\boldH} \chi_T( z) F_T^\Xi(z)dz \\
\nonumber&\le& \int_{\Gi\bk\boldH} \chi_T( z) F_T^{\G}(z)dz = \int_{\G\bk\boldH} \left(\sum_{\g\in\Gi\bk \G}\chi_T(\g z)\right)F_T^{\G}(z)dz \\
\label{FGTbound}&=&\|F_T^{\G}\|^2_{\G},
\eea
this time losing no information! Note that we have again used crucially the fact that $\Gi\subset \Xi\subset\G$.

As $\Psi_{\vep}^\Xi$ is an $\vep$-approximation to the identity on a two-dimensional space, we can choose it so that
\be\label{PsiGbound}
\|\Psi^\Xi_{\vep}\|_\Xi\ll {1\over\vep},
\ee
where the implied constant is independent of $\Xi$.
Combining \eqref{FGTbound} with \eqref{PsiGbound}, inserting into \eqref{almostThere}, and carrying out the same computation with $L_T$ replacing $K_T$, we have
\benn
\eqref{Errrr}\ll {1\over\vep} T^{s_1} \|F_3^{\G}\|_{\G} \ll {1\over\vep} T^{s_1},
\eenn
since $b<3$. 
The case $\gl_1\ge1/4$ is similar.

This completes the proof of the Lemma.
\end{proof}

The method of proof allows a much more general statement. Instead of pulling out only the contribution from the base eigenvalue $\gl_0$ in \eqref{Ginnerprod}, we can take more terms. Let $\gt\in(1/2,\gd)$ be a fixed constant, and denote the eigenvalues in $\Spec(\Xi\bk\boldH)$ below $\gt(1-\gt)$ by
\be\label{nameEvals}
\gd(1-\gd)=\gl_0<\gl_1\le \gl_2\le\dots\le\gl_M<\gt(1-\gt).
\ee
Corresponding to each point eigenvalue $\gl_j$ is a normalized $L^2$ eigenfunction, $\phi_j$. 
It may be the case that $\phi_j$ (which a priori is only $\Xi$-automorphic) is also automorphic with respect to  some group $\G_j$ which satisfies $\Xi\subset \G_j\subset\G$. Then $\phi_j$ is an ``oldform'' on $\Xi$ and $\phi_j\in L^2(\G_j\bk\boldH)$. This means $\gl_j=s_j(1-s_j)\in\Spec(\G_j\bk\boldH)$ and the same analysis as above gives
\benn
\hat{F}_T^\Xi(\gl_j)\hat{\Psi}^\Xi(\gl_j) ={ c_j T^{s_j}\over [\G_j :\Xi]} + O(T^{1/2}),
\eenn
where $c_j$ depends on $\G_j$ but not on $\Xi$. After extracting these lower order terms, the remaining error is simply $O({1\over\vep}T^{\gt})$. 
Suppressing the precise dependence on $\vep$ (which will be fixed for the remainder of this section),
we have proved

\begin{lemma}\label{moreterms}

Let $\Gi\subset \Xi\subset \G$ with $[\G:\Xi]<\infty$ and $w_T$ be defined by \eqref{wTdef}.  Fix $\gt\in(1/2,\gd)$ and let the eigenvalues in $\Spec(\Xi\bk\boldH)$ below $\gt(1-\gt)$ be denoted as in  \eqref{nameEvals}. 
For each $j=1,\dots,M$ let $\G_j$ denote a group satisfying $\Xi\subset \G_j\subset\G$ such that $\gl_j\in\Spec(\G_j\bk\boldH)$.

 Then for any fixed $\g\in\G$, there exists a constant $c_{\G}>0$ depending only on $\G$ and $\vep$, and constants $c_j$ depending on $\G_j$ and $\vep$, such that
\be\label{gBigSum}
\sum_{\xi\in\Gi\bk \Xi}w_T(\xi\g^{-1})={c_{\G} T^{\gd}\over [\G:\Xi] }+\sum_{j=1}^M{c_{j} T^{s_j}\over [\G_j:\Xi] }+ O_\vep\left(
T^{\gt}\right).
\ee

\end{lemma}

\subsection{Sums over $a_n$}

Recall the notation $f(\g)=c^2+d^2$ for a matrix $\g\in\G$ having bottom row $(c,d)$.
For $n\ge1$ let 
\be\label{aNdef}
a_n(T):=\sum_{{\g\in\Gi\bk\G}\atop{f(\g)=n}} w_T(\g)
\ee
be a smoothed count for the number of elements in our orbit having height bounded by $T$ and $f$-value exactly equal to $n$.

Recall that  
$\Gi\subset\G$
is the 
group which stabilizes $f(\g)=c^2+d^2$ in the sense that for $\g'\in\Gi$, we have $f(\g'\g)=f(\g)$.
Recall from \eqref{def:Gq} the principal ``congruence'' subgroup of level $q$ 
\benn
\G(q):=\{\g\in\G:\g\equiv I(q)\}.
\eenn
This group stabilizes all $\g$ $\mod  q$, i.e. if $\g'\in\G(q)$ then $\g'\g\equiv\g(q)$. Similarly, let
\be\label{G0qdef}
\G_1(q):=\{\g\in\G:\g\equiv\mattwo{ 1}{*}{0}{ 1}(\mod q)\}
\ee
be the subgroup of $\G$ which stabilizes $f(\mod  q )$, i.e. if $g\in\G_1(q)$ then $f(g\g)\equiv f(\g)(q)$. The 
inclusions $\Gi\subset\G_1(q)$
and $\G(q)\subset\G_1(q)\subset\G$
 are immediate. In particular, $\Spec(\G_1(q)\bk\boldH)\subset\Spec(\G(q)\bk\boldH)$ and $\G_1(q)$ inherits the spectral gap properties of $\G(q)$.

Fix $q\ge1$ square-free and consider
\be\label{ansSSum}
\sum_{n\equiv0(q)} a_n(T).
\ee

Insert \eqref{aNdef} into \eqref{ansSSum} and decompose $\g\in\Gi\bk\G$ into $\g=\xi \g_1$ with $\xi\in \Gi\bk \G_1(q)$ and $\g_1\in\G_1(q)\bk\G$:
\bea
\nonumber\sum_{n\equiv0(q)} a_n(T) &=& \sum_{{\g\in\Gi\bk\G}\atop{f(\g)\equiv0(q)}} w_T(\g) =\sum_{\g_1\in\G_1(q)\bk\G} \sum_{{\xi\in\Gi\bk\G_1(q)}\atop{f(\xi\g_1)\equiv0(q)}} w_T(\xi\g_1)\\
\label{gg0game}&=&  \sum_{{\g_1\in\G_1(q)\bk\G}\atop{\atop{f(\g_1)\equiv0(q)}}}\left( \sum_{{\xi\in\Gi\bk\G_1(q)}} w_T(\xi\g_1)\right),
\eea
since $f\mod q$ is invariant under $\G_1(q)$. Apply Lemma \ref{moreterms} with $\Xi=\G_1(q)$ to the inner sum in \eqref{gg0game} to prove

\begin{proposition}\label{anThm}
Let $T>1$, $q\ge1$ be square-free and $a_n(T)$ be defined by \eqref{aNdef}.  Let $\gt$ be the spectral gap of $\G$ and let the eigenvalues in $\Spec(\G_1(q)\bk\boldH)$ below $\gt(1-\gt)$ be denoted by
\benn
\gd(1-\gd)=\gl_0<\gl_1(q)\le \gl_2(q)\le\dots\le\gl_{M(q)}(q)<\gt(1-\gt).
\eenn

For each $j=1,\dots,M(q)$ let $q_j$ denote a divisor $q_j|q$ such that 
\[
\gl_j(q)=s_j(q)(1-s_j(q))\in\Spec(\G_1(q_j)\bk\boldH).
\]

 Then  there exists a constant $c_{\G}>0$ depending on $\G$ and $\vep$, and constants $c_j$ depending on $q_j$, such that
\be\label{aBigSum}
\sum_{n\equiv0(q)} a_n(T)  =|\Or_q|\left({c_{\G} T^{\gd}\over [\G:\G_1(q)] }+\sum_{j=1}^{M(q)}{c_{j} T^{s_j(q)}\over [\G_1(q_j):\G_1(q)] }+ O\left({1\over\vep}
T^{\gt}\right)\right),
\ee
where
\be\label{Oq0def}
|\Or_q| =|\{(c,d)(\mod q)\in\Or\mid c^2+d^2\equiv0(q)\}|= \sum_{{\g_1\in\G_1(q)\bk\G}\atop{f(\g_1)\equiv0(q)}}1.
\ee
\end{proposition}

\section{Proof of the Main Theorem}\label{finish}

In this section we prove Theorem \ref{mainThm}.

\subsection{Part (1)}

We begin by proving \eqref{mainthm1}. Let 
\[
H_{\vep}(T):= \sum_{n}a_n(T),
\]
with $a_n(T)$ defined in \eqref{aNdef} and implicitly dependent on $\vep$ via $w_{T,\vep}$ in \eqref{wTdef}. 
Recall the definition of $\Or(T)$ from \eqref{def:OT}. By \eqref{dubyabound} we have
\[
|\Or(T/(1+\vep))|\le H_{\vep}(T)\le |\Or(T/(1-\vep))|,
\]
or equivalently,
\be\label{thisworksQ}
H_{\vep}(T (1-\vep))\le|\Or(T)|\le H_{\vep}(T (1+\vep)).
\ee
Apply 
Lemma \ref{GG0Game}
 with $\Xi=\G$ and $\g=I$:
\be\label{thatworks}
H_{\vep}(T) = (c_{\G}+O(\vep))T^{\gd} + O({1\over\vep}T^{s_1}),
\ee
where $\gl_1=s_1(1-s_1)$ is the first eigenvalue above the base in $\Spec(\GbkH)$, satisfying $\gd>s_1$. Combining \eqref{thisworksQ} with \eqref{thatworks} we elementarily arrive at \eqref{mainthm1} by an appropriate choice of $\vep$.

\subsection{Part (2)}

Our only remaining task is to verify all of the conditions necessary to apply a combinatorial sieve to $\A=\{a_n(T)\}$ (see  \S \ref{combsieve}).
By \eqref{dubyabound}, the sequence $a_n(T)$ defined by \eqref{aNdef} satisfies \eqref{aNwithTbound} with $T(1+\vep)$ replacing $T$. Throughout the rest of the section, $\vep$ will be a fixed small constant. Anyway since the bound on $R$ in \eqref{Rbound} depends on $\log T$, this difference is irrelevant.

Let $\gt\in[1/2,\gd)$ be a spectral gap for $\G$. There are at most finitely many primes for which Strong Approximation fails, and also finitely many primes at which the corresponding spectrum fails to have a $\gt$-gap. Let $\B$ be the product of the primes in these finite ``ramified'' sets. 

For ease of exposition, assume first that $q\ge 1$ is square-free and relatively prime to $\B$. 
Then the projection
\be\label{ontoG}
\G\longrightarrow \SL_2(\Z/q\Z)\text{ is onto,}
\ee
and
\be\label{specAgain}
\Spec(\G(q)\bk\boldH)\cap (0,\gt(1-\gt)) = \Spec(\GbkH)\cap (0,\gt(1-\gt)).
\ee
Apply Proposition \ref{anThm} to \eqref{Aqs}, and infer from \eqref{specAgain} that we can take $q_j=1$ for all $j$:
\benn
|\A_q|=|\Or_q|\left( {1\over [\G:\G_1(q)]}\left(c_{\G} T^{\gd}+\sum_{j}c_j T^{s_j} \right) + O_{\vep}(T^{\gt})  \right),
\eenn
where all $c_j$'s only depend on $\G$ (and $\vep$, which is fixed) but not on $q$. 
 Then $|\A_q|=\gw(q)\X+r(q)$, with
\bea
\nonumber
\gw(q)&=&{|\Or_q|\over [\G:\G_1(q)]},\\
\nonumber
\X&=&c_{\G}T^{\gd}+\sum_j c_j T^{s_j},
\\
\nonumber
& \text{\hskip-3.7in and }&\\
\label{rQbound}
r(q)&=&O({|\Or_q| T^{\gt}}).
\eea
Trivially $\gw(1)=1$. By \eqref{Oq0def} and \eqref{ontoG} it is clear that $\gw$ is multiplicative, and for $q=p$ a prime we can compute $\gw(p)$ exactly. 

As $p$ is unfamified, $|\Or_p|$ counts the number of $(c,d)\mod p$ with $(c,d)\neq(0,0)$ and $c^2+d^2\equiv0(p)$. The last equation is equivalent to $c^2\equiv-d^2(p)$, which has no solutions if $-1$ is  not a square mod $p$, that is, when  $p\equiv 3(4)$. In the opposite case, this cardinality is easily computed by hand. Similarly, $[\G:\G_1(p)]$  counts the number of $(c,d)\mod p$ with $(c,d)\neq(0,0)$. Thus we have
\bea
\label{OrpEval}
|\Or_p| &=&
\threecase{1}{if $p=2$}{2(p-1)}{if $p\equiv1(4)$}{0}{if $p\equiv3(4)$,} \\
\nonumber
[\G:\G_1(p)]& =& p^2-1, 
\\
\nonumber
& \text{\hskip-2.4in and so }&\\
\nonumber
\gw(p) &=&
{|\Or_p|
\over
[\G:\G_1(p)]}
=
\threecase{1/3}{if $p=2$}{2/(p+1)}{if $p\equiv1(4)$}{0}{if $p\equiv3(4)$.}
\eea
Then the sieve condition (S\ref{gwsat}) is obvious, and   (S\ref{localDensBound}) 
follows from
$$
\prod_{p<z\atop p\equiv1(4)} \left(1-{2\over p}\right)^{-1}\sim \gk \log z
,
$$
a classical exercise (see e.g. \cite{Landau1953}).
Inserting \eqref{OrpEval} into  \eqref{rQbound} gives
\be\label{rest}
r(q)=O(q T^{\gt}),
\ee
and since $\X\asymp T^{\gd}$, (S\ref{rqErr}) requires
\benn
\sum_{q\le Q} |r(q)| = O (Q^2 T^{\gt}) \ll T^{\gd/(1+\gep)},
\eenn
for any $\gep>0$. This is satisfied for
\be\label{Qfinal}
Q= T^{(\gd-\gt)\over 2(1+\gep)}.
\ee
Then inputting \eqref{Qfinal} into  \eqref{Rbound}, together with \eqref{eq:Sz}, gives \eqref{Part2}
for 
\benn
R>2 \log T/\log Q = 4 / (\gd-\gt),
\eenn
since $\gep>0$ was arbitrary.

This completes the analysis of the affine linear sieve in the case $q$ is ``unramified''.

\subsection{Ramified places}

Let $q=q'q''\ge1$ be square-free, with $q'|\B$
and $(q'',\B)=1$. 
Notice that $\B$ being the product of a  finite number of primes means there are only finitely many possible values of $q'$.

Then by Theorem \ref{strongApp} the projection of $\G$ in $\SL_2(\Z/q\Z)$ is 
\be\label{ontoG2}
G_{q'}\times \SL_2(\Z/q''\Z), 
\ee
where $G_{q'}$ is the projection of $\G$ in $\SL_2(\Z/q'\Z)$.

By Theorem \ref{BGSspec}, if $\gl<\gt(1-\gt)$ is in $\Spec(\G(q)\bk\boldH)
$ but not in $\Spec(\GbkH)$, then $\gl\in\Spec(\G(q')\bk\boldH)$. 

Then applying Proposition \ref{anThm} and using the fact that 
\[
[\G:\Xi]=[\G:\G_j][\G_j:\Xi]
\] 
for $\Xi\subset \G_j\subset\G$, we have
\benn
|\A_q|=|\Or_q|\left( {1\over [\G:\G_1(q)]}\left(\X + \X_q\right) + O_{\vep}(T^{\gt})  \right),
\eenn
where 
\[
\X =
c_{\G} T^{\gd} +
\sum_{{\gl_j\in\Spec(\GbkH)-\{\gl_0\}}\atop{\gl_j=s_j(1-s_j)<\gt(1-\gt)}}    c_j T^{s_j}, 
\]
and
\[
\X_{q'} =  [\G:\G_1(q')] 
\sum_{{\gl_j^*\in\Spec(\G_1(q')\bk\boldH)_{new}}  \atop{\gl_j^*=s_j^*(1-s_j^*)<\gt(1-\gt)}    }                c_j^* T^{s_j^*}.
\]
Here $c_j$'s depend on $\G$ and not on $q$, while $c_j^*$'s depend on $q'$ but not on $q''$.
Clearly $\X_{q'}\ll\X^{1-\eta}$ for some $\eta>0$.
By \eqref{ontoG2}, $|\Or_q|=|\Or_{q'}||\Or_{q''}|$ and $[\G_1(q''):\G_1(q)] = [\G:\G_1(q')]$, which is independent of $q''$.

 Therefore we have the expression
\[
|\A_q| = \gw(q)(\X+\X_{q'}) + r(q),
\]
where 
\[
\gw(q) ={|\Or_{q}|\over [\G:\G_1(q)] } = {|\Or_{q'}| \over [\G:\G_1(q')]}{|\Or_{q''}|\over [\G:\G_1(q'')] } = \gw(q')\gw(q''),
\]
and $r(q)$ satisfies \eqref{rest}. The rest of the analysis follows 
as before, completing the proof of the Main Theorem.

\appendix

\section{Proof of Theorem  \ref{sieve}}\label{AxSieve}

In this appendix we derive Theorem \ref{sieve} from the more standard sieve setting.
%
See 
\cite{Iwaniec1996,IwaniecKowalski2004}.

As before, let $\A=\{a_n\}$ be our sequence of nonnegative numbers with $a_n=0$ for $n$ exceeding a parameter $T$, and let $\B$ be the product of a finite set of ``bad'' primes.
Recall that $\X$ is an approximation to $|\A|$: 
\benn
\X\asymp\sum_n a_n \asymp T^{\gd}
\eenn
and 
that for $q=q'q''$ square-free with $q'|\B$ and $(q'',\B)=1$ we have
\[
|\A_q|=\sum_{n\equiv0(q)} a_n = \gw(q) (\X + \X_{q'})+r(q).
\]
It is assumed that for any divisor $q'|\B$, the factor $\X_{q'}$ is a power less than the main term $\X$. Also $\gw$ is multiplicative away from $\B$, i.e. $\gw(q)=\gw(q')\prod_{p|q''}\gw(p)$.
 Let $z\ge2$ be a parameter (a small power of $\X$) 
 and define 
\beann
P(z)&:=&\prod_{{p<z}}p. 
\eeann
Assume $z$ is large enough (by taking $T$ large enough) so that $\B|P(z)$.
Consider a sum of the form
\beann
S(z)&:=&\sum_{(n,P(z))=1} a_n.
\eeann

Notice that if  $z=T^{\ga}$, $n\le T$ and $(n,P(z))=1$ then $n$ has at most $1/\ga$ prime factors. In this way, $S(z)$ counts the number of  $1/\ga$-almost primes. 

By M\"{o}bius inversion, we have
\beann
S(z) &=& \sum_n a_n \sum_{q|(n,P(z))} \mu(q)\\
 &=& \sum_{q|P(z)} \mu(q) \sum_{n\equiv0(q)} a_n\\
 &=& \sum_{q|P(z)} \mu(q) \left(  \gw(q) (\X + \X_{q'})+r(q)      \right) \\
 &=& \gS_1+\gS_2+\gS_3,
 \eeann
with 
\[
\gS_1 = \sum_{q|P(z)} \mu(q) \gw(q) \X=\X\prod_{{p
<z
}\atop{p\nmid\B}}(1-\gw(p))\times\sum_{q'|\B}\mu(q')\gw(q'),
\]
\[
\gS_2 = \sum_{q|P(z)} \mu(q) \gw(q) \X_{q'}= \prod_{{p<z}\atop{p\nmid\B}}(1-\gw(p))\times \left(\sum_{q'|\B}\mu(q')\gw(q')\X_{q'}\right),
\]
and
\[
\gS_3 = \sum_{q|P(z)} \mu(q) r(q).
\]

 Let 
 \[
 S^*(z)= \gS_1+\gS_3.
 \]

\begin{theorem}[See Iwaniec \cite{Iwaniec1996}, Theorem 3.7 on page 63]\label{thm:sieve}
Let $Q>e^{2K}$ where $K$ is the constant appearing in \eqref{localDensBound} and define
\beann
V(z)&=& \prod_{{p\le z}\atop{p\nmid\B}}(1-\gw(p)),\\
R(Q) &=& \sum_{q\le Q}|r(q)|\\
f(s) &=& 2e^{\g}\log(s-1)/s\text{, for $s\in[2,4]$},\\
F(s) &=& 2e^{\g}/s\text{, for $s\in[1,3]$, and }\\
D&=&c K^{11} (\log\log\log Q)^3(\log\log Q)^{-1}.
\eeann
Here $\g=.577\dots$ is the Euler constant and $c=3.591\dots$ solves $(c/e)^c=e$.
Then 
\benn
(f(s)-D)\X V(z)-R(Q) \le S^*(z)\le (F(s)+D)\X V(z)+R(Q),
\eenn
where $s=\log Q/\log z$.
\end{theorem}

The choice
$
Q=T^{{\gd-\gt\over 2(1+\gep)} }
$
gives
\benn
D\ll (\log\log\log T)^3(\log\log T)^{-1}\to0,
\eenn
so we  can take 
\be\label{eq:s}
s=2(1+\gep), 
\ee
giving $f(s)>D$ for $T$ sufficiently large.

Then with $z=Q^{1/s}=T^{{\gd-\gt\over 4(1+\gep)^2} }$ we have 
\benn
S^*(z)\asymp {\X\over\log\X} \asymp T^{\gd}/\log T\asymp S(z),
\eenn
since $\gS_2\ll\X^{1-\eta}$.
As $\gep>0$ is arbitrary, this captures $R$-almost primes with  any $R>\delthe$.

\begin{remark}\label{rmkA2} The set $Z_R$ is genuinely the set of integers having at most $R$ prime factors, not just $R$ prime factors outside of $\B$.
\end{remark}

\section{Proof of Lemma \ref{lem:ODE}}\label{technic}

We require the following simple lemma from the theory of inhomogeneous ODEs, in particular the method of  variation of parameters.

\begin{lemma}\label{lem:phi}
Let $\gl\ge0$ and suppose $f$ and $g$ are functions satisfying 
\benn
-y^2\frac{\partial^2}{\partial y^2}f(y) - \gl f(y) = g(y).
\eenn

Assume $\gl=s(1-s)\neq1/4$. Then there exist constants $\ga$ and $\gb$ such that
\benn
f(y) = \ga y^s + \gb y^{1-s} + u(y) y^s + v(y) y^{1-s},
\eenn
where
\be\label{udef}
u(y) = (1-2s)^{-1} \int_{1/T}^y w^{-1-s} g(w)dw
\ee
and
\benn
v(y) = (2s-1)^{-1} \int_{1/T}^y w^{s-2} g(w)dw.
\eenn
If $\gl=1/4$ then 
\benn
f(y) = \ga y^{1/2} + \gb y^{1/2}\log y + u(y) y^{1/2} + v(y) y^{1/2}\log y,
\eenn
where
\beann
u(y)&=&\int_{1/T}^y w^{-3/2}\log(w) g(w) dw\text{, and}\\
v(y)&=&-\int_{1/T}^y w^{-3/2} g(w) dw.
\eeann
\end{lemma}

\begin{proof}
Elementary calculus. For details see \cite{MyThesis}.
\end{proof}

Fix any $\phi\in L^2(\GbkH)$ and any $\gl\ge0$. For simplicity assume $\gl\neq1/4$ (the calculation in the opposite case is similar). Consider  the left hand side of \eqref{ODE} and unfold:
\beann
\< F_T, \phi \> &=& \int_{\GbkH} F_T(z)\phi(z) dz \\
 &=& \int_{1/T}^{\infty} \left(\int_N\bar{\phi}(z)dx\right)\frac{dy}{y^2} ,
\eeann
where again $N$ is an interval corresponding to the width  of the cusp at infinity.

Let 
\benn
f(y) = \int_N\bar{\phi}(z) dx,
\eenn
so that 
\be\label{littlefInt}
\<F_T,\phi\>  = \int_{1/T}^{\infty} f(y) \frac{dy}{y^2}.
\ee
Let $g$ be defined by:
\be\label{gdef}
g(y):= -y^2 \frac{\partial^2}{\partial y^2} f(y) - \gl f(y) =  \int_N(\gD-\gl)\phi(z) dx .
\ee
By Lemma \ref{lem:phi},
\benn
f(y)= \ga y^s + \gb y^{1-s} + y^su(y) + y^{1-s}v(y).
\eenn
The first two terms are the homogenous solution and the last two are the perturbation.
Of course inserting the homogenous component of $f$ into \eqref{littlefInt} we have the main term in \eqref{ODE}
\benn
\int_{1/T}^{\infty} \left( \ga y^s +\gb y^{1-s}\right) \frac{dy}{y^2} = A T^s + B T^{1-s}.
\eenn
Thus it remains to show that $I, II \ll_{\gl,T} \|(\gD-\gl)\phi\|$, where
\benn
I=\int_{1/T}^{\infty} y^s u(y)\frac{dy}{y^2}, \text{ and } II=\int_{1/T}^{\infty} y^{1-s} v(y)\frac{dy}{y^2}
\eenn
are the contributions from the perturbation.
Integrate $I$ by parts and recall from \eqref{udef} that $u(1/T)=0$:
\be\label{eq:I}
I=\left. u(y)\frac{y^{s-1}}{s-1} \right|_{y\to\infty} - \int_{1/T}^{\infty} \frac{y^{s-1}}{s-1}\left(\frac{y^{-1-s}g(y)}{1-2s}\right)dy.
\ee
Modulo constants,  insert \eqref{gdef} into the last integral of \eqref{eq:I} and apply Cauchy-Shwarz:
\bea
\nonumber 
& & \int_{1/T}^{\infty} \int_N \left(\frac{(\gD-\gl)\phi(z)}{y}\right)\left(\frac{1}{y}\right)dx\,dy\\ 
\nonumber 
 & \ll & \left( \int_{1/T}^{\infty} \int_N \left| \frac{(\gD-\gl)\phi(z)}{y} \right|^2 dx\,dy \right )^{1/2} \left( \int_{1/T}^{\infty} \int_N \left| \frac{1}{y} \right|^2 dx\,dy\right )^{1/2}  \\ 
\label{CSobst}  & \ll_{\gl,T,\G} & \|(\gD-\gl)\phi\|\, \left(|N| \,T\right)^{1/2}\\
\nonumber   & \ll_{\gl,T,\G} & \|(\gD-\gl)\phi\|,
 \eea
where $|N|$ is the length of $N$.
Here
we used the fact that the box 
$$
N\times[1/T,\infty]\subset\boldH
$$ 
is contained 
in a union of finitely many (depending on $T$) fundamental 
domains for $\G\backslash\boldH$. 
\begin{remark}\label{baad}
The appearance of $|N|$ in \eqref{CSobst} (via the use of Cauchy-Shwarz)  is the {\bf most severe} obstruction to removing the assumption that $\G$ stabilize infinity. If there is no cusp at infinity, then $|N|=\infty$ and our analysis fails. See Remark \ref{whybaad}.
\end{remark}

For the first part of \eqref{eq:I}, we need the following bound:
\beann
u(y) & = & \frac{1}{1-2s}\int_{1/T}^{y} w^{-1-s}g(w)dw 
\\
 & = &
  \frac{1}{1-2s}\int_{1/T}^{y} \int_N \frac{(\gD-\gl)\phi(x+\i w)}{w}\cdot w^{-s}dx\,dw 
 \\
 & \ll_{\gl} & 
 \left( \int_{1/T}^{y} \int_N \left| \frac{(\gD-\gl)\phi(x+\i w)}{w} \right|^2 dx\,dw \right)^{1/2} 
 \\
 &&
 \times\left( \int_{1/T}^{y} \int_N \left|  w^{-s} \right|^2 dx\,dw \right)^{1/2}    
 \\
 & \ll_{\gl,T} & \|(\gD-\gl)\phi\| \twocase{}{y^{1-2s}+T^{2s-1}}{if $s>1/2$}{\log y +\log T}{if $\Re(s)=1/2$.} 
\eeann
So
\benn
\lim_{y\to\infty} u(y) y^{s-1} \ll_{\gl,T} \|(\gD-\gl)\phi\|.
\eenn

The integral $II$ is handled identically and we are done.

\subsection*{Note added in proof:} In joint  work with Hee Oh \cite{KontorovichOh2008}, we circumvent the Main Identity
to prove the Main Theorem without the assumption that $\Gi$ is nontrivial.
 Instead of the Main Identity, we prove the equidistribution of long horocycle flows on the unit tangent bundle of an infinite-volume Riemann surface of constant negative
curvature, and then use this equidistribution to count. Furthermore, we replace the $\beta$ sieve by the weighted sieve of Diamond-Halberstam-Richert \cite{DiamondHalberstamRichert1988}, which gives better numbers under nearly identical hypotheses, and execute the sieve for various other choices of $f$ in \eqref{fIs}. We also use these methods in \cite{KontorovichOh2009} to count the number of circles in an Apollonian packing of bounded curvature, and discuss various Diophantine properties of integral 
Apollonian
packings.

\bibliographystyle{alpha}

\bibliography{Submit}

\end{document}